\documentclass[11pt]{amsart}
\usepackage{latexsym,amsmath,amsthm,amsfonts}
\usepackage[psamsfonts]{amssymb}

\newtheorem{theorem}{\sc Theorem}[section]
\newtheorem{corollary}[theorem]{\sc Corollary}
\newtheorem{lemma}[theorem]{\sc Lemma}

\theoremstyle{definition}

\theoremstyle{remark}

\numberwithin{equation}{section}

\def\shat{{\mathbb S}}
\def\wh{\widehat}
\def\eps{\varepsilon}
\def\C{\mbox{$\mathcal C$}}
\def\X{\mbox{$\mathcal X$}}
\def\RR{\mathbb R}
\def\EE{\mathbb{E}_{f}}
\def\PP{\mathbb{P}}

\def\argmin{\mathop{\rm arg\, min}}

\def\dim{\mathop{\rm dim}}
\def\pen{\mathop{\rm pen}}
\def\egal{\stackrel{{\rm def}}{=}}

\begin{document}
\baselineskip=18pt

\title[Aggregation for regression learning]{Aggregation for regression
learning}
\author[Bunea, Tsybakov, and Wegkamp]{Florentina Bunea$^{\dag 1}$}
\address{Florentina Bunea\\ Department of Statistics,   Florida State
University,   Tallahassee,   Florida.}
\email{\rm bunea@stat.fsu.edu}
\author[]{Alexandre B. Tsybakov}
\address{Alexandre B. Tsybakov\\ Laboratoire de Probabilit\'es et
Mod\`eles Al\'eatoires, Universit\'e Paris VI, France.}
\email{\rm tsybakov@ccr.jussieu.fr}
\author[]{Marten H. Wegkamp$^\dag$}
\address{Marten H. Wegkamp\\ Department of Statistics,   Florida State
University,   Tallahassee,   Florida.}
\email{\rm wegkamp@stat.fsu.edu}

\thanks{$^1$ Corresponding author}
\thanks{$^\dag$Research partially supported by NSF grant DMS 0406049 }

\subjclass{Primary 62G08, Secondary 62C20, 62G05, 62G20}

\keywords{aggregation, minimax risk, model selection,  nonparametric regression, oracle inequalities, penalized least squares, statistical
learning}

\date{October  2004}

\bibliographystyle{amsplain}

\begin{abstract}
This paper studies statistical aggregation procedures in regression
setting. A motivating factor is the existence of many different
methods of estimation, leading to possibly competing estimators.

We consider here three different types of aggregation: model
selection (MS) aggregation, convex (C) aggregation and  linear (L)
aggregation. The objective of (MS) is  to select the optimal single
estimator from the list; that of (C) is to select  the optimal
convex combination of the given estimators; and that of (L) is to
select the optimal
 linear combination of the given estimators. We are
interested in evaluating the rates of convergence of the excess
risks of the estimators obtained by these procedures. Our approach
is motivated by recent minimax results in  Nemirovski (2000) and
Tsybakov (2003).

There exist competing aggregation procedures achieving optimal
convergence separately for each one of (MS), (C) and (L) cases.
Since the bounds in these results are not directly comparable with
each other, we suggest an alternative solution. We prove that all
the three optimal bounds can be nearly achieved via a single
``universal" aggregation procedure. We propose such a procedure
which consists in mixing of the initial estimators with the weights
obtained by penalized least squares. Two different penalities are
considered: one of them is related to hard thresholding techniques,
the second one is a data dependent $L_1$-type penalty.

\end{abstract}

\maketitle

\section{Introduction}

\noindent In this paper we study aggregation procedures and their performance for regression models. Let $\mathcal{D}_n = \{
(X_1,Y_1),\ldots,(X_n,Y_n) \}$ be a sample of independent random pairs $(X_i,Y_i)$ with
\begin{equation}
\label{eq:model} Y_i = f(X_i) + W_i, \quad i=1,\ldots,n,
\end{equation}
where $f:{\mathcal X} \to \RR$ is an unknown regression function to be
estimated, ${\mathcal X}$ is a Borel subset of $\RR^d$, the $X_i$'s are
either random vectors with probability measure $\mu$ supported on
${\mathcal X}$ or fixed elements in $\X$, and the errors $W_i$ are zero mean
random variables, conditionally on 
the $X_i$'s.

Aggregation of arbitrary estimators in regression models has
recently received increasing attention: Nemirovski (2000), Juditsky
and Nemirovski (2000), Yang (2000, 2001, 2004), Catoni (2001),
Gy\"orfi {\em et al.} (2002), Wegkamp (2003), Tsybakov (2003),
Birg\'e (2003). A motivating factor is the existence of many
different methods of estimation, leading to possibly competing
estimators. Local polynomial kernel smoothing
 methods and penalized least squares or likelihood estimators
(which include B-splines and wavelet type estimators) are two
classes of methods that cover the major trends in nonparametric
estimation in regression. When no method is a clear winner, one may
prefer to combine different estimators obtained via different
methods. Furthermore, within each method one can obtain competing
estimators for different values of the smoothing parameter (the
bandwidth in kernel procedures and, for the other examples, the
calibrating constant in the penalty term or, correspondingly, the
threshold value). This is usually the case when adaptive estimation
is considered. In all these situations we are faced with a large
collection of concurrent estimators $\wh{f}_1,\ldots, \wh{f}_{M}$. A
natural idea is then to look for a new, improved, estimator
$\widetilde{f}$ constructed by combining $\wh{f}_1,\ldots,
\wh{f}_{M}$ in a suitable way. Such an estimator $\widetilde{f}$ is
called {\it aggregate} and its construction is called aggregation.

There exist three main aggregation problems: model selection (MS)
aggregation, convex (C) aggregation and linear (L) aggregation. They
are discussed in detail by Nemirovski (2000). The objective of (MS)
is  to select the optimal (in a sense to be defined) single
estimator from the list; that of (C) is to select  the optimal
convex combination of the given estimators; and that of (L) is to
select the optimal linear combination of the given estimators.

In this paper we consider a more general setup for the (MS), (C) and (L)
aggregation problems, following Tsybakov (2003).
Namely, we do not
restrict aggregates to be of the form of model selectors, convex or
linear combinations of the original estimators. Instead, we only require
that aggregates should be estimators that mimic the model selection,
convex or linear oracles.  This allows us to construct more powerful
aggregates. To give precise definitions, denote by $\|g\| = \left(\int
g^2(x)\mu(dx)\right)^{1/2}$ the norm of a function $g$ in $L_2(\RR^d,
\mu)$ and set ${\sf f}_{\lambda} = \sum_{j = 1}^{M} \lambda_j\wh{f}_j$ for
any $\lambda=(\lambda_1,\dots,\lambda_M)\in \RR^M$. The performance of an
aggregate $\widetilde{f}$ used to estimate a function $f \in L_2(\RR^d,
\mu)$ can be judged against the following mathematical target:
\begin{eqnarray}\label{ner}
\mathbb E_f \|\widetilde{f} - f\|^{2}\leq \inf_{\lambda \in H^M} \mathbb
E_f\|{\sf f}_{\lambda} - f\|^{2} + \Delta_{n,
M},
\end{eqnarray}
where $\Delta_{n, M}\ge 0$ is a remainder term {\it independent of $f$}
characterizing the price to pay for aggregation, and the set $H^M$ is
either the whole $\RR^M$ (for linear aggregation),  or the simplex
$\Lambda^M=\left\{\lambda=(\lambda_1,\dots,\lambda_M) \in \mathbb R^M:
\lambda_j
\ge 0,\ \sum_{j = 1}^{M} \lambda_j \leq 1\right\}$
(for convex aggregation), or the set of $M$ vertices of $\Lambda^M$ (for
model selection aggregation). Here and later $\mathbb E_f$ denotes the
expectation with respect to the joint distribution of
$(X_1,Y_1),\ldots,(X_n,Y_n)$ under model (\ref{eq:model}). The random
functions
${\sf f}_{\lambda}$ attaining  $\inf_{\lambda \in H^M} \mathbb
E_f\|{\sf f}_{\lambda} - f\|^{2}$
 in (\ref{ner}) for the three values taken by $H^M$ are called (L), (C)
and (MS) oracles, respectively. Note that these minimizers are not
estimators since they depend on the true $f$.

We say that the aggregate $\widetilde{f}$ mimics the (L), (C) or (MS)
oracle if it satisfies (\ref{ner}) for the corresponding set $H^M$, with
the minimal possible price for aggregation $\Delta_{n, M}$. Minimal
possible values $\Delta_{n, M}$ for the three problems can be defined via a
minimax setting and they are called optimal rates of aggregation [Tsybakov (2003)]
and further denoted by $\psi_{n,M}$.
As shown in Tsybakov (2003),
for the Gaussian regression model we have, under mild conditions
\begin{equation}
\label{oragg}
\psi_{n, M} \asymp \left\{\begin{array}{cl}
{M}/{n}   & \quad \mbox{for  (L) aggregation},\\\\
{M}/{n}   & \quad \mbox{for  (C) aggregation,  if }  M\le \sqrt{n},\\\\
\sqrt{ \left\{ \log(1+M/\sqrt{n})\right\}/n} & \quad
\mbox{for  (C) aggregation,  if }  M  >  \sqrt{n},\\\\
({\log M})/{n}   &  \quad \mbox{for  (MS) aggregation}.
\end{array} \right.
\end{equation}
This implies that linear aggregation has the highest price, (MS)
aggregation has the lowest one, and convex aggregation occupies an
intermediate place. The oracle risks on the right in (\ref{ner})
satisfy a reversed inequality:
\[
\inf_{1 \leq j \leq M} \mathbb E_f\|{f}_{j} - f\|^{2}  \geq \inf_{\lambda
\in \Lambda^M} \mathbb E_f\|{\sf f}_{\lambda} - f\|^{2} \geq
    \inf_{\lambda \in
\mathbb R^M} \mathbb E_f\|{\sf f}_{\lambda} - f \|^{2},
\]
since the sets over which the infima are taken are nested. Thus, the
bound (\ref{ner}) for (MS) aggregation realizes the trade-off between the
largest oracle risk and the smallest remainder term. The bound
(\ref{ner}) for (L) aggregation realizes the trade-off between the smallest
oracle risk and the largest remainder term. The bound (\ref{ner}) for (C)
aggregation realizes the trade-off between an intermediate oracle risk
and intermediate remainder term. If the number of estimators to be
aggregated is small, $M \leq \sqrt{n}$, the  remainder term in the (C) bound
is identical to that in the (L) bound, but the  oracle risk in the (L)
bound is always superior to that in the (C) bound.
 Thus (L) aggregation is preferable to
(C) aggregation in this case, but no comparison can be made with
(MS) aggregation. If the number of estimators to be aggregated is
large, $M > \sqrt n$, the remainder term in the (L) bound becomes
too large, but, in a strict sense, there is no winner among the
three aggregation techniques. The question how to choose the best
among them remains open.

The ideal oracle inequality (\ref{ner}) is available only for some
special cases. See Catoni (2001) for (MS) aggregation in Gaussian
regression; Nemirovski (2000), Juditsky and Nemirovski (2000),
Tsybakov (2003) for (C) aggregation with $M>\sqrt{n}$; and Tsybakov
(2003) for (L) aggregation with known marginal measure $\mu$ and for
(C) aggregation with $M\le\sqrt{n}$.  For more general situations
there exist less precise results of the type
\begin{equation}
\label{ner1}
 \mathbb E_f \|\widetilde{f} - f\|^{2}
 \leq C_0\inf_{\lambda \in H^M} \mathbb E_f\|{\sf f}_{\lambda} - f\|^{2} +
\Delta_{n, M},
\end{equation}
where $C_0 > 1$ is a constant independent of $f$ and $n$, and
$\Delta_{n, M}$ is a remainder term, not necessarily having the same
behavior in $n$ and $M$ as the optimal one $\psi_{n,M}$. A
disadvantage of (\ref{ner1}) over (\ref{ner}) is that, when the
oracle risk $R^* = \inf_{\lambda \in H^M} \mathbb E_f\|{\sf
f}_{\lambda} - f\|^{2}$ is large, the additional term $(C_0-1)R^*$
on the right-hand side of (\ref{ner1}) may be much larger than the
remainder term $ \Delta_{n, M}$, thus substantially spoiling the
convergence properties. This effect is less pronounced if $C_0 = 1
+\eps $ for some arbitrarily small $\eps>0$ or for $\eps=\eps_n\to0$
as $n\to\infty$.

Bounds of the type (\ref{ner1}) in regression problems have been
obtained by many authors mainly for the model selection case (when
$H^M$ is the set of vertices of the simplex $\Lambda^M$), see, for
example, Kneip (1994), Barron {\em et al.} (1999), Lugosi and Nobel
(1999), Catoni (2001), Gy\"{o}rfi {\em et al.} (2002), Baraud (2000,
2002), Bartlett {\em et al.} (2002), Wegkamp (2003), Birg\'e (2003),
Bunea (2004), Bunea and Wegkamp (2004), and the references cited in
these works. Most of the papers on model selection treat particular
restricted families of estimators, such as orthogonal series
estimators, spline estimators, etc. An interesting recent
development due to Leung and Barron (2004) covers model selection
for all estimators admitting Stein's unbiased estimation of the
risk. There are relatively few results on (MS) aggregation when the
estimators are allowed to be arbitrary, see Catoni (2001), Yang
(2000, 2001, 2002), Gy\"{o}rfi {\em et al.} (2002), Wegkamp (2003),
Birg\'e (2003), and Tsybakov (2003). Here we make the standard
assumption that $\wh{f}_1,\ldots, \wh{f}_{M}$ are uniformly bounded,
but otherwise they can be arbitrary.

Various convex aggregation procedures for nonparametric regression
have emerged in the last decade. They include bootstrap based
methods, as suggested by LeBlanc and Tibshirani (1996)
 and cross-validation based stacking, as in Wolpert (1992) or Breiman (1996).
 The literature on oracle inequalities of the type (\ref{ner}) and
(\ref{ner1}) for the (C) aggregation case is not nearly as large as
the one on model selection. Juditsky and Nemirovski (2000),
Nemirovski (2000) propose a stochastic approximation algorithm that
achieves the bound (\ref{ner})  for (C) aggregation with optimal
rate $\psi_{n,M}$ in the case $M>n/\log n$. They also show that the
bound is achieved by usual (non-penalized) least squares convex
aggregation. Yang (2000, 2001, 2004) suggest several methods of
convex aggregation, in particular ARM (adaptive regression by
mixing). He proves bounds of the form (\ref{ner1}) with constants
$C_0$ that are typically much larger than 1 and with rates
$\Delta_{n,M}$ that can be equal or approximately equal to the
optimal rates $\psi_{n,M}$ when $M$ is  a power of $n$.  Audibert
(2003) establishes (\ref{ner}) for a PAC-Bayesian method of convex
aggregation with almost optimal rates, up to a logarithmic factor.
Birg\'e (2003) suggests a convex aggregation method satisfying
(\ref{ner1}) with a constant $C_0$ that can be much greater than 1
and with a rate that is optimal for $M>\sqrt{n}$ and suboptimal for
$M\le\sqrt{n}$. On the other hand, Koltchinskii (2004, Section 8)
proves (\ref{ner}) for a convex aggregate $\widetilde{f}$ with
optimal rate for $M\le\sqrt{n}$ and with almost optimal rate for
$M>\sqrt{n}$.

Linear aggregation procedures have received substantially less
attention. For regression models with random design, a procedure
achieving the bound (\ref{ner}) with optimal rate $\psi_{n,M}$ of
(L) aggregation can be found in Tsybakov (2003). For Gaussian white
noise models, linear aggregation has been discussed earlier by
Nemirovski (2000).

Aggregation procedures are typically based on sample splitting. The initial sample ${\mathcal D}_n$ is divided into two independent subsamples
${\mathcal D}_m^1$ and ${\mathcal D}_\ell^2$ of sizes $m$ and $\ell$, respectively, where $m\gg \ell$ and $m+\ell=n$. The first subsample
${\mathcal D}_m^1$ (called training sample) is used to construct estimators $\wh{f}_{1}, \dots, \wh{f}_{M}$ and the second subsample ${\mathcal
D}_\ell^2$ (called learning sample) is used to aggregate them ({\em i.e.}, to construct ${\widetilde f}$). In this paper we do not consider
sample splitting schemes but rather deal with an idealized scheme. Following Nemirovski (2000),
 the first subsample is fixed and thus
instead of estimators $\wh{f}_{1}, \dots, \wh{f}_{M}$, we have fixed
functions $f_1, \dots, f_M$. That is, we focus our attention on
learning. Our aim is to find estimators based on the sample
${\mathcal D}_n$ that would mimic simultaneously the linear, convex
and model selection oracles with the fastest possible rates (or,
equivalently, with the smallest possible remainder terms
$\Delta_{n,M}$). A passage to the initial model is straightforward:
it is enough to condition on the first subsample, to use the
learning bounds of the type (\ref{ner}), (\ref{ner1}) obtained for
the idealized scheme, and then to take expectations of both sides of
the inequalities over the distribution of the whole sample
${\mathcal D}_n$.

Another interpretation of aggregation of fixed functions $f_1, \dots,
f_M$ is related to parametric regression for linear models of dimension
$M$, where $M$ can be very large or increasing with $n$. In fact, assume
that both $X_i$ and  $\wh{f}_j = f_j$ are fixed (non-random), and
consider the linear regression model with design matrix $\left(
f_{j}(X_i) \right)_{1\le i\le n,\ 1\le j\le M}$ and the empirical
counterpart of
the norm $\|\cdot\|$ defined by
\[ \| f\|_n = \left( \frac1n \sum_{i=1}^n f^2(X_i) \right)^{1/2}. \]
Then, for $H^M = \Lambda^M$ or $H^M = \RR^M$, the value $\inf_{\lambda
\in H^M} \mathbb \|{\sf f}_{\lambda} - f\|_n^{2}$ represents the best least
squares approximation of an unknown function $f$ at points $X_i$ by the
convex or linear span, respectively, of the columns of the design
matrix. Consequently, estimators $\widetilde{f}$ satisfying oracle
inequalities of the form
\begin{equation}
\label{ner2}
 \mathbb E_f \|\widetilde{f} - f\|_n^{2}
 \leq C_0\inf_{\lambda \in H^M} \mathbb \|{\sf f}_{\lambda} - f\|^{2}_n +
\Delta_{n, M}
\end{equation}
mimic the best linear/convex least-squares approximation of $f$ in a
parametric regression framework, provided $C_0\ge 1$ is close to
$1$. In (\ref{ner2}), $\Delta_{n, M}$ can be interpreted as the
price to pay for the dimension $M$ of the regression model, and we
will show that (for an appropriate choice of the aggregate
$\widetilde{f}$) $\Delta_{n, M}=\psi_{n,M}$, where $\psi_{n,M}$ is
the optimal rate of aggregation as defined in (\ref{oragg}). For the
case of linear aggregation, this can be viewed in the spirit of
earlier work on linear models with growing dimension $M$ [Yohai and
Maronna (1979), Portnoy (1984)], but here we obtain non-asymptotic
results and our risk is defined in terms of the regression functions
and not in terms of their parameters.

Given the existence of competing aggregation procedures achieving
either optimal  (MS), or (C), or (L) bounds, there is an ongoing
discussion as to which procedure is the best one. Since this cannot
be decided by merely comparing the optimal bounds, we suggest an
alternative solution. We show that all the three optimal  (MS), (C)
and (L) bounds can be nearly achieved via a single aggregation
procedure. Consequently, the smallest of the three will be achieved.
Our answer will thus meet the desiderata of both model selection and
model averaging.

The procedures that we suggest for aggregation are based on penalized
least squares.
We consider two penalties that can be associated
with soft thresholding ($L_1$ or Lasso type penalty) and with hard
thresholding, respectively.

In Section 3.1 we show that a hard threshold
aggregate satisfies
inequalities of the type (\ref{ner2}), with $C_0$
arbitrarily close to 1, and with the optimal  remainder term $\psi_{n, M}$.
We establish the oracle inequalities for all three sets $H^M$ under consideration,
hence showing that the hard threshold
aggregate achieves simultaneously the (MS), (C) and (L) bounds
when the empirical norm $\|\cdot\|_n$ is used to define the risk.

In Section 3.2 we study the performance of a slightly different hard
threshold aggregate under the $L_2(\RR^d, \mu)$ norm. We show that
this aggregate satisfies simultaneously the oracle inequalities of
the type (\ref{ner1}) corresponding to the (MS) and (C) bounds, with
a remainder term $\Delta_{n, M}$ that possibly differs from the
optimal $\psi_{n, M}$ in a logarithmic factor,  and with $C_0$
arbitrarily close to 1.


In Section 4 we study aggregation with the $L_1$ penalty and we obtain (\ref{ner2}) simultaneously for the (MS), (C) and (L) cases,  with $C_0$
arbitrarily close to 1 and with a remainder term $\Delta_{n, M}$ that differs from the optimal $\psi_{n, M}$ only in a logarithmic factor.

Finally, we study lower bounds for (MS) and (L) aggregation in the fixed design case in Section 5, complementing the results obtained for the
random design case by Tsybakov (2003).\\

\bigskip

\section{Notation and assumptions}

\noindent The following two assumptions on the regression model (\ref{eq:model}) are
supposed to be satisfied throughout the paper.\\

\noindent{\sc Assumption (A1)} {\it The random variables $W_i$ are
independent and Gaussian} $N(0,\sigma^2)$.\\

\noindent {\sc Assumption (A2)} {\it The functions $f:{\mathcal X}\to\RR$
and $f_j:{\mathcal X}\to\RR, \ j=1,\dots, M$, with
 $M\ge 2$, belong to the class ${\mathcal F}_0$ of uniformly bounded functions defined by
 \[ {\mathcal F}_0 \egal \Big\{g:{\mathcal X}\to\RR\, \Big| \
 \sup_{x\in\X} |g(x)|\le L  \Big\}
\]
where $L < \infty$ is a constant that is not necessarily known to the statistician.}\\

The normality assumption (A1) on the distribution of errors is convenient since
we need certain exponential tail bounds in the proofs (see Lemma
\ref{LM} below). For example, bounded regression can be easily
incorporated in this framework using maximal inequalities due to
Talagrand (1994a, b) and
Panchenko (2003). More generally,   subgaussian errors are allowed at the cost of increasing
technicalities, see  Van de Geer (2000).
In order to retain a transparent presentation of both the results and
proofs, we confine ourselves to the Gaussian regression framework.\\

For any $\lambda= (\lambda_1,\dots,\lambda_M)\in \RR^M$, define
$$ {\sf f}_{\lambda}(x) = \sum_{j=1}^M \lambda_j f_j(x).$$
The functions $f_j$ can be viewed as estimators of $f$ constructed from a training sample (see the Introduction). Here we consider the ideal
situation in which they are fixed, {\em i.e.}, we concentrate on learning only. The learning method that we propose is based on aggregating  the
$f_j$'s via penalized least squares.

For each $\lambda=(\lambda_1,\ldots,\lambda_M)\in\RR^M$, let $M(\lambda)$
denote the number of non-zero coordinates of $\lambda$:
\[
M(\lambda)= \sum_{j=1}^M I\{ \lambda_j\ne 0\} = {\rm Card} \ J(\lambda)
\]
where $I\{ \cdot\}$ denotes the indicator function, and $J(\lambda)=\left
\{j\in\{1,\ldots,M\}:\ \lambda_j \ne 0\right\}$. Introduce the
residual sum of squares
\[ \wh{S} (\lambda) = \frac1n \sum_{i=1}^n \{ Y_i - {\sf f}_{\lambda}(X_i) \}^2.\]
Given a penalty term pen($\lambda$), the penalized least squares
estimator $\wh{\lambda}
= (\wh{\lambda}_1,\dots,\wh{\lambda}_M)$ is defined by
\begin{eqnarray}
\label{lse}
\wh{\lambda}
= \arg \min_{\lambda\in\RR^M} \left\{ \wh{S}(\lambda) +
\text{pen}(\lambda) \right\},
\end{eqnarray}
 which renders in turn the aggregated
estimator
$$
\widetilde{f}(x)={\sf f}_{\wh{\lambda}}(x).
$$
Since the vector $\wh{\lambda}$ can take any values in $\RR^M$, the
aggregate $\widetilde{f}$ is not a model selector in the traditional sense, nor is
it  necessarily a convex combination of the functions $f_j$.
Nevertheless, we will show that it mimics the  (MS), (C) and (L) oracles
when one
of the following two penalties is used:
\begin{eqnarray}
\label{pen1} \text{pen}(\lambda) = K_1 \frac{ M({\lambda})}{n} \log\left(
1 + \frac{M}{M({\lambda})\vee 1}\right)
\end{eqnarray}
or
\begin{eqnarray}
\label{pen2} \text{pen}(\lambda) =
\sum_{j=1}^M r_{n,j}|\lambda_j|,
\end{eqnarray}
where $K_1>0$ is a constant independent of $M, n$, and $r_{n,j}$'s are the
data-dependent weights defined in (\ref{pen5}).\\

We refer to the penalty in (\ref{pen1}) as  {\it hard threshold penalty}. This is motivated by the well known fact that, in the sequence space
model ({\em i.e.}, when the functions $f_1,\dots,f_M$ are orthonormal with respect to the scalar product induced by the norm $\|\cdot\|_n$), the
penalty $\pen(\lambda) \sim M(\lambda)$ leads to $\wh{\lambda}_j$'s that are hard thresholded values of the $Y_j$'s (see, for instance,
 H\"ardle {\em et al.} (1998), page 138).
Our penalty (\ref{pen1}) is not exactly of that form, but it differs from
it only in a logarithmic factor.

The penalty (\ref{pen2}), again in the sequence space model,
leads to $\wh{\lambda}_j$'s that are soft thresholded values of
$Y_j$'s. We will call it therefore  {\it soft threshold penalty} or
$L_1$-{\it penalty}.
Penalized least squares estimators with soft threshold
penalty $\text{pen}(\lambda) \sim \sum_{j=1}^M |\lambda_j|$ are
closely related to
Lasso-type  estimators [Tibshirani (1996), Efron {\em et al.} (2004)].
Our results show that, with $r_{n,j}$'s defined by (\ref{pen5}),
the soft threshold penalty allows near optimal
aggregation. The same is true for the hard threshold penalty (\ref{pen1})
under
somewhat different conditions. \\

In what follows, we denote by $C, C_1, C_2, \dots$ finite positive
constants, possibly different on different occasions.
\bigskip

\section{Near optimal aggregation with the hard threshold penalty}
\subsection{\sc The fixed design case}
In this section we show that the penalized least squares estimator using
a penalty of the form (\ref{pen1}) achieves simultaneously the (MS),
(L), and (C) bounds of the form (\ref{ner2}) with the correct rates
$\Delta_{n,M}=\psi_{n,M}$. Consequently, the smallest bound is achieved by our aggregate. The
results of this section are established for the empirical loss $\|
\widetilde{f} - f \|_n^2$. The next theorem presents an oracle inequality which
implies all the three bounds.\\

\begin{theorem}\label{theorem:prima}
Let $X_i \in {\mathcal X}$, $i=1,\dots,n$, be fixed.
Let $\widetilde{f}$ be the penalized least squares estimate defined in
(\ref{lse}) with penalty
(\ref{pen1}).
There exist constants $C_1,C_2>0$  such that for all $a>1$, for
$K_{1}=K_0 a\sigma^2$, with $K_0>0$ large enough, and for all integers
$n\ge 1$ and $M\ge 2$,
\begin{eqnarray}
\label{oracle} && \EE \| \widetilde{f} - f \|_n^2 \\ &&\le
\inf_{\lambda\in\RR^M} \left\{ \frac{a+1}{a-1}\| {\sf f}_{\lambda} - f \|_n^2
 + \ C_1 a \sigma^2 \frac{ M({\lambda})}{n} \log\left( 1 +
\frac{M}{M({\lambda})\vee 1}\right) \right\}
 +  C_2\frac{  a \sigma^2}{n}. \nonumber
\end{eqnarray}
\end{theorem}

\medskip

This theorem is proved in Section 3.3. The following three corollaries
present bounds of the form (\ref{ner2}) for (MS), (L), and (C)
aggregation, respectively.

\medskip

\begin{corollary}[MS]\label{cor:MS}
Let the assumptions of Theorem \ref{theorem:prima} be satisfied.
Then there exists a constant $C_3>0$ such that for all $\eps>0$, for
$K_{1}=K_{1}(\eps,\sigma^2)$ large enough and for all integers $n\ge
1$ and $M\ge 2$,
\begin{eqnarray}
\EE \| \widetilde{f} - f \|_n^2   \le (1+\eps)\inf_{1\le j\le M}  \|
f_j - f \|_n^2 + C_3 \sigma^2 \left(1+{\eps}^{-1}\right) \frac{ \log
M}{n}. \nonumber
\end{eqnarray}
\end{corollary}
\begin{proof}
Since the infimum on the right of (\ref{oracle}) is taken over all
$\lambda\in\RR^M$, the bound easily follows by  considering only the
subset consisting of the $M$ vertices
$(\lambda_1,\ldots,\lambda_M)=(1,0,\ldots,0)$,
$(0,1,0,\ldots,0),\ldots,(0,\ldots,0,1)$ in $\Lambda^M$, and by
putting $a=1+2/\eps$.
\end{proof}

\bigskip

\begin{corollary}[L]
\label{corL} Let the assumptions of Theorem \ref{theorem:prima} be
satisfied. Then there exists a constant $C_3>0$ such that for all
$\eps>0$, for $K_{1}=K_{1}(\eps,\sigma^2)$ large enough and for all
integers $n\ge 1$ and $M\ge 2$,
\begin{eqnarray}
\EE \| \widetilde{f} - f \|_n^2   \le (1+\eps)\inf_{\lambda\in\RR^M}
\| {\sf f}_{\lambda} - f  \|_n^2 + C_3 \sigma^2
\left(1+{\eps}^{-1}\right) \frac{ M}{n}. \nonumber
\end{eqnarray}
\end{corollary}
\begin{proof}
Since $x\mapsto x \log (1+M /x)$ is increasing for $1\le x\le M$,
\begin{eqnarray*}
\label{L1} \sup_{\lambda\in\RR^M}
 \frac{ M({\lambda})}{n} \log\left( 1 + \frac{M}{M({\lambda})\vee
1}\right) = \frac{M}{n} \log2.
\end{eqnarray*} The result then follows from (\ref{oracle}) with $a=1+2/\eps$.
\end{proof}

\bigskip

\begin{corollary}[C]\label{cor:C}
Let the assumptions of Theorem \ref{theorem:prima} be satisfied.
Then there exists a constant $C_3'>0$ depending on $L$ and
$\sigma^2$ such that for all $\eps>0$, for
$K_{1}=K_{1}(\eps,\sigma^2)$ large enough and for all integers $n\ge
1$ and $M\ge 2$,
\begin{eqnarray}
\EE \| \widetilde{f} - f \|_n^2 \le
(1+\eps)\inf_{\lambda\in\Lambda^M} \| {\sf f}_{\lambda} - f \|_n^2 +
C_3'\left(1+\eps+{\eps}^{-1}\right) \psi_n^C(M), \nonumber
\end{eqnarray}
where
\begin{eqnarray*} \psi_n^C(M) = \begin{cases} M/n & \text{ if }
M\le\sqrt{n},\\
\sqrt{\{ \log( 1+ M/\sqrt{n}) \} / n}   & \text{ if } M> \sqrt{n}.
\end{cases}\end{eqnarray*}
\end{corollary}
\begin{proof}
For $M\le \sqrt{n}$ the result follows from Corollary \ref{corL}.
Assume now that $M> \sqrt{n}$ and let $m$ be the integer part of
$$x_{n,M} = \sqrt{\frac{n\log 2}{\log( 1+ M/\sqrt{n}) }}.$$
Clearly, $0\le m\le x_{n,M} \le M$. First, consider the case
$m\ge1$.
 Denote by $\C$ the set of functions $h$ of the form
\[  h(x)= \frac1m \sum_{j=1}^M {k_j} f_j(x), \ k_j\in\{0,1,\ldots,m\}, \ \sum_{j=1}^m k_j \le m .
\]
The following approximation result  can be obtained by the
``Maurey argument" (see, for example, Barron (1993), Lemma 1, or
Nemirovski (2000), pages 192, 193):
\begin{eqnarray}
\label{one} \min_{g\in\C} \| g-f\|_n^2 \le \min_{\lambda\in
\Lambda^M} \| {\sf f}_{\lambda}- f\|_n^2 + \frac{L^2}{m}.
\end{eqnarray}
For completeness, we give the proof of (\ref{one}) in the
Appendix. Since $M(\lambda)\le m\le x_{n,M}$ for the vectors
$\lambda$ corresponding to $g\in\C$, and since $x\mapsto x\log
\left(1+\frac{M}{x}\right)$ is increasing for $1\le x\le M$, we
get from (\ref{oracle}):
\[ \EE \|\widetilde{f}-f\|_n^2 \le \inf_{g\in\C }
\left\{ \frac{a+1}{a-1} \| g-f\|_n^2 + C_1 a\sigma^2 \
\frac{x_{n,M}}{n} \log \left( 1+ \frac{M}{x_{n,M}} \right)\right\}
+ \frac{C_2a\sigma^2}{n}.
\]
Using this inequality, (\ref{one}) and the fact that $m=\lfloor
x_{n,M}\rfloor \ge x_{n,M}/2$ for $x_{n,M}\ge1$, we obtain
\begin{eqnarray}
\label{two} \EE\|\widetilde{f}-f\|_n^2 &\le& \frac{a+1}{a-1}
\inf_{\lambda\in\Lambda^M}\| {\sf f}_{\lambda}-f\|_n^2 +
\left(\frac{a+1}{a-1}\right)\frac{2L^2}{x_{n,M}}
\\
&&+ \ C_1 a\sigma^2 \ \frac{x_{n,M}}{n} \log\left(
1+\frac{M}{x_{n,M}}\right) + \frac{C_2a\sigma^2}{n}.\nonumber
\end{eqnarray}
We use this bound for all choices of $\lambda\in\Lambda^M$ with
$m\ge M(\lambda)\ne0$. For  $m=0$, we only need to consider the
singular case $\lambda=0$ as $M(\lambda)=0$ if and only if
$\lambda=0$. Note that for $m=0$, we have $1/x_{n,M} \ge 1$, and
we use the trivial upper bound
$$
\frac{a+1}{a-1} \| f\|_n^2  + \frac{C_2a\sigma^2}{n} \le
\left(\frac{a+1}{a-1} L^2  + C_2a\sigma^2\right)\left(\frac{\log(
1+ M/\sqrt{n})}{n\log 2}\right)^{1/2}
$$
for the right-hand side of (\ref{oracle}).

To complete the proof of the Corollary, it remains to put
$a=1+2/\eps$ and to note that
$$
\log\left( 1+\frac{M}{x_{n,M}}\right)\le 2\log
\left(1+\frac{M}{\sqrt{n}}\right),
$$
in view of the elementary inequality $\log\Big(1+(\log 2)^{-1/2} y
\sqrt{\log(1+y)} \Big) \le 2\log(1+y)$, for all $y\ge 1$.
\end{proof}

\medskip

\noindent We remark now that the aggregate considered in Theorem
\ref{theorem:prima} satisfies also the bounds ``in probability" that
are similar in spirit to (\ref{oracle}) and its corollaries.

\begin{theorem}\label{theorem:inprob1}
Let $X_i \in {\mathcal X}$, $i=1,\dots,n$, be fixed.
Let $\widetilde{f}$ be the penalized least squares estimate defined in
(\ref{lse}) with penalty
(\ref{pen1}).
There exist constants $C_1,L_1, L_2>0$  such that for all $a>1$, for
$K_{1}=K_0 a\sigma^2$, with $K_0>0$ large enough, and for all integers
$n\ge 1$, $M\ge 2$ and any $\delta > 0$,
\begin{eqnarray}
\label{oracle1} && \PP \left (  \| \widetilde{f} - f \|_n^2  \geq
\inf_{\lambda\in\RR^M} \left\{ \frac{a+1}{a-1}\| {\sf f}_{\lambda} - f \|_n^2
 + \ C_1 a \sigma^2 \frac{ M({\lambda})}{n} \log\left( 1 +
\frac{M}{M({\lambda})\vee 1}\right) \right\}
 +  \delta \right ) \\
&&\hspace{3cm} \leq  L_1 \exp\left(- L_2\frac{n\delta}{a\sigma^2}\right). \nonumber
\end{eqnarray}
\end{theorem}

As in the case of Theorem \ref{theorem:prima}, we can consequently
obtain the analogues of Corollaries \ref{cor:MS} - \ref{cor:C}, by
replacing the infimum in  (\ref{oracle1}) by its particular form for
the cases (MS), (L) and (C), respectively. We do not include each
case, for brevity.

\bigskip

\subsection{\sc The random design case}

In  this subsection we show that an oracle inequality similar to
(\ref{oracle}) continues to hold if  the empirical norm $\|\cdot\|_n$ is
replaced
by the $L_2(\RR^d,\mu)$ norm $\|\cdot\|$. This result is more difficult
to obtain and we do not achieve exactly the same bounds.

We need to restrict minimization of the penalized sum of squares
to a bounded set in $\RR^M$. Define, for any $T>0$,
\[ \Lambda_{M,T} =\left\{ \lambda\in\RR^M:\ \sum_{j=1}^M |\lambda_j|\le
T\right\}.\]
 The penalty term needs to be chosen slightly larger than before:
\begin{eqnarray}
\label{pen4} \text{pen}(\lambda) = K_{1} \frac{ M({\lambda})}{n}
\log\left( 1 + \frac{M\vee n}{M({\lambda})\vee 1}\right)
\end{eqnarray}
for some large $K_1>0$. We note that here $K_1$ is not necessarily
the same as in (\ref{pen1}), we just use the same notation for
factors in the penalty term.

\begin{theorem}\label{th2}
Assume that $X_1,\ldots,X_n$ are independent random variables with
common probability measure $\mu$. Let $T<\infty$ be fixed, and set
$$B= L^2(T+1)^2.$$ Let $\widetilde{f}=f_{\wh{\lambda}}$ where
\[
\wh{\lambda} = \argmin_{\lambda\in\Lambda_{M,T}} \{\wh{S}(\lambda)+
{\rm pen}(\lambda)\}
\]
with the penalty given in (\ref{pen4}). Then there exist constants
$C_1,C_2>0$  such that for all $a>1$, for $K_{1}=K_{1}(a, B,
\sigma^2)$ large enough, and for all integers $n\ge 1$ and $M\ge 2$,
\begin{eqnarray}
\label{rades} && \EE \| \widetilde{f} - f \|^2\\ && \le
\inf_{\lambda\in\Lambda_{M,T}} \left\{ \frac{a+1}{a-1}\| {\sf
f}_{\lambda} - f \|^2 + \ C_1 a \sigma^2 \frac{ M({\lambda})}{n}
\log\left( 1 + \frac{M\vee n}{M({\lambda})\vee 1}\right) \right\}
 +  C_2\frac{a(\sigma^2+B)}{n}. \nonumber
\end{eqnarray}
\end{theorem}

\medskip

Because of the slight increase in the penalty, the remainder term in
(\ref{rades}) is somewhat larger than the one given in
(\ref{oracle}): we now have $M\vee n$ in place of $M$ under the
logarithm.

As corollaries, one obtains the following (MS)  and (C) bounds for
the estimator $\widetilde{f}$ defined in Theorem \ref{th2}.\\
\begin{corollary}[MS]
\label{cor:MS2}
Let the assumptions of Theorem \ref{th2} be satisfied and $T\ge 1$.
Then there exists a constant $C>0$ such that for all $\eps>0$, for
$K_{1}=K_{1}(\eps,\sigma^2)$ large enough and for all integers $n\ge
1$ and $M\ge 2$,
\begin{eqnarray}
\EE \| \widetilde{f} - f \|^2   \le (1+\eps)\inf_{1\le j\le M}  \|
f_j - f \|^2 + C \sigma^2 \left(1+{\eps}^{-1}\right) \frac{ \log
(M\vee n)}{n}. \nonumber
\end{eqnarray}
\end{corollary}

\bigskip

\begin{corollary}[C]\label{cor:C2}
Let the assumptions of Theorem \ref{th2} be satisfied and $T\ge1$.
Then there exists a constant $C'>0$ depending on $L$ and $\sigma^2$
such that for all $\eps>0$, for $K_{1}=K_{1}(\eps,\sigma^2)$ large
enough and for all integers $n\ge 1$ and $M\ge 2$,
\begin{eqnarray}
\EE \| \widetilde{f} - f \|^2 \le (1+\eps)\inf_{\lambda\in\Lambda^M}
\| {\sf f}_{\lambda} - f \|^2 + C' \left(1+\eps + {\eps}^{-1}\right)
\widetilde{\psi}_n^C(M), \nonumber
\end{eqnarray}
where
\begin{eqnarray*} \widetilde{\psi}_n^C(M) = \begin{cases} (M \log n) / n
& \text{ if }
M\le\sqrt{n},\\
\sqrt{\{ \log( 1+ (M\vee n)/\sqrt{n}) \}/ n}   & \text{ if } M>
\sqrt{n}.
\end{cases}\end{eqnarray*}
\end{corollary}

\bigskip

As compared to Corollaries \ref{cor:MS} and \ref{cor:C}, these
results present slightly different rates of convergence: here the
factor $\log M$ is replaced by $\log n$ for values $M<n$. The proofs
are omitted since Corollaries \ref{cor:MS2} and \ref{cor:C2} readily
follow from the oracle inequality (\ref{rades}) and the fact that
$\Lambda^{M}\subset \Lambda_{M,T}$ for $T\ge 1$ via an argument
similar to the proofs of Corollaries \ref{cor:MS} and \ref{cor:C}.

\bigskip

\subsection{\sc Proof of Theorem \ref{theorem:prima}}
Let $\lambda$ be a fixed, but arbitrary point in $\RR^M$. Define for all
$1\le m \le M$,
$$
A_m(\lambda) = \{ \bar\lambda = \lambda^\prime - \lambda \in\RR^M:\
M(\lambda^\prime)=m\} .
$$
Let $J_k, \ k=1,\dots,\binom{M}{m}$, be all the subsets of
$\{1,\dots,M\}$ of cardinality $m$. Define
$$ A_{m,k}(\lambda) =\left \{ \bar\lambda =
(\bar\lambda_1,\dots,\bar\lambda_M)\in A_m(\lambda):\
\lambda^\prime_j \ne 0 \Leftrightarrow j\in  J_k \right\}
$$
where $\lambda^\prime_j = \bar\lambda_j + \lambda_j$. The collection
$\left \{ A_{m,k}(\lambda):\ 1 \le k\le \binom{M}{m} \right\}$ forms a
partition of the set $A_m(\lambda)$.
 Furthermore, define affine  subspaces of $\RR^n$ of the form
\[
B_{m,k}(\lambda)= \left\{ h =\left( {\sf f}_{\bar\lambda}(X_1), \dots,
{\sf f}_{\bar\lambda}(X_n) \right) \in \RR^n:\ \bar\lambda \in A_{m,k}(\lambda)
\right\}
\]
and let $\Pi_{m,k}^{\lambda}W$ denote the projection of the vector
$W=(W_1,\ldots,W_n)$ onto $B_{m,k}(\lambda)$. Clearly, $\dim
(B_{m,k}(\lambda)) \le m$. Finally, we define for each $\gamma \in\RR^M$,
\begin{eqnarray*} V_n(\gamma) =
 \frac1n \sum_{i=1}^n W_i \frac{{\sf f}_{\gamma}(X_i)}{\| {\sf f}_{\gamma}\|_n} \
\text{ if }\
\| {\sf f}_{\gamma}\|_n\ne 0,
\end{eqnarray*}
and $V_n(\gamma)\egal0$, otherwise.
\medskip
\begin{lemma}\label{een} For all $a>1, b>0$ and $\lambda\in\RR^M$, we have
\begin{eqnarray*}
\| \widetilde{f} - f  \|_n^2 & \le & \frac{1+b}{b}\frac{a}{a-1} \| {\sf f}_{\lambda} - f
\|_n^2 + \frac{a}{a-1} K_{1} \frac{ M({\lambda})}{n} \log\left( 1 +
\frac{M}{M({\lambda})\vee 1}\right)
\nonumber\\
&&+ \ \frac{a}{a-1}\max_{1\le m \le M} \max_{1\le k\le \binom{M}{m}} \left\{
 (a+b) \| \Pi_{m,k}^{\lambda}W \|_n^2 - \frac{K_{1} m}{n}  \log\left(1+
\frac{M}{m\vee 1}\right) \right\} \\ && +\ \frac{a(a+b)}{a-1} V_n^2(\lambda).
\end{eqnarray*}
\end{lemma}
\begin{proof}
By the definition of $\widehat \lambda$, for any $\lambda\in\RR^M
$,
\[ \wh{S}(\wh{\lambda})+ \text{pen}(\wh{\lambda}) \le \wh{S}(\lambda)+
\text{pen}(\lambda).\]
Rewriting this inequality yields
\begin{eqnarray}\nonumber
\| \widetilde{f} -f \|_n^2 &\le& \| {\sf f}_{\lambda} - f \|_n^2 + 2\left< W, \widetilde{f}-
{\sf f}_{\lambda} \right>_n + \text{pen}(\lambda) -
\text{pen}(\wh{\lambda}),\end{eqnarray} where $<\cdot,\cdot>_n$ denotes
the scalar product associated with the norm $\|\cdot\|_n$. Since $
\|\widetilde{f}-{\sf f}_{\lambda}\|_n=0$ implies that $\left< W, \widetilde{f}- {\sf f}_{\lambda}
\right>_n=0$,  we find
\begin{eqnarray*}
\| \widetilde{f} -f \|_n^2
 &\le& \| {\sf f}_{\lambda} - f \|_n^2 + 2V_n(\wh{\lambda}-\lambda)
\| \widetilde{f}- {\sf f}_{\lambda} \|_n  + \text{pen}(\lambda)
- \text{pen}(\wh{\lambda})\\
&\le& \| {\sf f}_{\lambda} - f \|_n^2 + 2V_n(\wh{\lambda}-\lambda) \| \widetilde{f}- f \|_n
+ 2V_n(\wh{\lambda}-\lambda) \|  {\sf f}_{\lambda} -f \|_n + \text{pen}(\lambda)
- \text{pen}(\wh{\lambda})\\
&\le& (1+\frac{1}{b}) \| {\sf f}_{\lambda} - f \|_n^2 + a
V_n^2(\wh{\lambda}-\lambda) + \frac1a \| \widetilde{f}- f \|_n^2 +
bV_n^2(\wh{\lambda}-\lambda)  +
\text{pen}(\lambda) - \text{pen}(\wh{\lambda}),
\end{eqnarray*}
where $a,b>0$ are arbitrary, and we used the inequality $2xy \le c x^2 +
y^2/c$ valid for all $x,y\in\RR$ and $c>0$. Consequently, for any
$a>1,b>0$, we find
\begin{eqnarray*}
\| \widetilde{f} -f \|_n^2 &\le& \frac{1+b}{b}\frac{a}{a-1} \| {\sf f}_{\lambda} - f
\|_n^2 + \frac{a}{a-1} \text{pen}(\lambda)\\ &&+ \frac{a}{a-1}(a+b)
V_n^2(\wh{\lambda}-\lambda) - \frac{a}{a-1} \text{pen}(\wh{\lambda}).
\end{eqnarray*}
Next, since $\RR^M = \bigcup_{m=0}^M \bigcup_{k=1}^{\binom{M}{m}}
A_{m,k}(\lambda)$, we find that
\begin{eqnarray*}
&& (a+b) V_n^2(\wh{\lambda}-\lambda)    -
\text{pen}(\wh{\lambda})\\
&&= (a+b) V_n^2(\wh{\lambda}-\lambda)
-  \text{pen}(\wh{\lambda}-\lambda+\lambda)\\
 &&\le
 \max_{0 \le m \le M} \max_{1\le k\le \binom{M}{m}} \max_{\bar\lambda
\in A_{m,k}(\lambda)  }\left\{ (a+b) V_n^2(\bar\lambda)   -
\text{pen}(\bar\lambda+\lambda) \right\}.
\end{eqnarray*}  It
remains to bound the term on the right in view of the last two displays.
The case $m=0$ is degenerate as
 $A_0(\lambda)=A_{0,1}(\lambda)=\{-\lambda\}$. Note that
 for $\bar\lambda=-\lambda$,
\[  (a+b) V_n^2(\bar\lambda)  -
\text{pen}(\bar\lambda+\lambda)
 = (a+b) V_n^2(\lambda) ,\] since $\pen(0)=0$ and
${\sf f}_{-\lambda}=-{\sf f}_{\lambda}$.
  For each $m\ge1$, we have
\begin{eqnarray*}
&&\max_{1\le k\le \binom{M}{m}} \max_{\bar\lambda \in A_{m,k}(\lambda)
}\left\{ (a+b) V_n^2(\bar\lambda)  - \text{pen}(\bar\lambda+\lambda)
\right\}
\\
&& \le \max_{1\le k\le \binom{M}{m}} \max_{\bar\lambda \in
A_{m,k}(\lambda) }\left\{
 (a+b) \| \Pi_{m,k}^{\lambda}W \|_n^2 - \text{pen}(\bar\lambda+\lambda)
\right\}\\
 &&\quad \text{ by the orthogonality of $W-\Pi_{m,k}^{\lambda}W$ and $
\left({\sf f}_{\bar\lambda}(X_1),\ldots, {\sf f}_{\bar\lambda}(X_n)
\right)$ for all $\bar\lambda\in A_{m,k}(\lambda)$}\\
&& = \max_{1\le k\le \binom{M}{m}} \left\{
 (a+b) \| \Pi_{m,k}^{\lambda}W \|_n^2 - \frac{K_{1}}{n} m \log\left(1 +
\frac{M}{m\vee 1}\right)
  \right\}
  \\
&&\quad\text{ in view of (\ref{pen4}) and since $M(\bar\lambda
+\lambda)=m$  for all
$\bar\lambda\in A_{m,k}(\lambda)$}.
\end{eqnarray*}
This concludes the proof of the lemma. \end{proof}

From now on, we take $a=b>1$. Since, by Assumption (A1), the
errors $W_i$ are normal $N(0,\sigma^2)$,  the standardized
statistic $n \sigma^{-2} \| \Pi_{m,k}^{\lambda}W \|_n^2$ has a
$\chi^2$ distribution with $m$ degrees of freedom for all $1\le
k\le \binom{M}{m}$. The following tail bound for such a statistic
will be useful.
\begin{lemma}\label{LM} Let $Z_d$ denote a random variable having the
$\chi^2$ distribution with $d$ degrees of freedom.
Then for all $x>0$,
\begin{eqnarray}
\PP\{ Z_d-d \ge x \sqrt{2d}  \} \le \exp\left(-\frac{x^2}{
2(1+x\sqrt{2/d})}\right).
\end{eqnarray}
\end{lemma}
\begin{proof}
See Cavalier {\em et al.} (2002),
equation (27) at page 857.
\end{proof}
\begin{lemma}
\label{twee}
There exists $C>0$  such
that, for any integer $n\ge 1$ and any $a>1$, $K_{1}=K_0a\sigma^2$ with
$K_0>0$ large enough,
\begin{eqnarray}
&& \label{3.10.1}
 \EE \max_{1\le m \le M} \max_{1\le k\le \binom{M}{m}} \left\{
2a \| \Pi_{m,k}^{\lambda}W \|_n^2 - \frac{K_{1}}{n} m \log\left( 1+
\frac{M}{m\vee 1}\right)
  \right\} \le C \frac{a\sigma^2}{ n },
\\
&& \label{3.10.2} \EE V_n^2(\lambda) \le  \frac{\sigma^2}{ n }.
  \end{eqnarray}
 \end{lemma}
 \begin{proof}
Inequality (\ref{3.10.2}) is trivial and we will prove only
(\ref{3.10.1}). For any $\delta>0$ we have
\begin{eqnarray*}
&& p_\delta \egal \PP\left[ \max_{1\le m \le M} \max_{1\le k\le
\binom{M}{m}} \left\{
 2a \| \Pi_{m,k}^{\lambda}W \|_n^2 - \frac{K_{1}}{n} m
\log\left(1+\frac{M}{m\vee 1}\right)
  \right\} \ge \delta \right] \\
 && \le \sum_{m=1}^M \sum_{k=1}^{\binom{M}{m}} \PP\left[
 2a \| \Pi_{m,k}^{\lambda}W \|_n^2 - \frac{K_{1}}{n} m \log\left( 1 +
\frac{M}{m\vee 1}\right) \ge \delta \right] \\
     && = \sum_{m=1}^M \sum_{k=1}^{\binom{M}{m}} \PP\left[
 Z_m \ge   \frac{K_{1}}{2a\sigma^2} {m }  \log\left(1+\frac{M}{m}\right)+
\frac{n \delta}{2a\sigma^2} \right]\\
 && = \sum_{m=1}^M
\binom{M}{m} \PP \left[ \frac{Z_m - m}{\sqrt{2m}} \ge
\frac{K_{1}}{2a\sigma^2} \frac{\sqrt{m}}{\sqrt{2} }
\log\left(1+\frac{M}{m}\right)
-\frac{\sqrt{m}}{\sqrt{2}} + \frac{n \delta}{2a\sigma^2 \sqrt{2m}}
 \right]\\
&& \le
\sum_{m=1}^M \binom{M}{m} \exp\left(- C_0 \left\{ \frac{ m K_{1} }{a\sigma^2}
\log\left(1+\frac{M}{m}\right) + \frac{n\delta}{a\sigma^2}
\right\} \right)
\end{eqnarray*}
by Lemma \ref{LM} for $K_{1}=K_0a \sigma^2$ with $K_0>0$ large enough and
some universal constant $C_0>0$. Using the crude bound $\binom{M}{m}  \le
\left({e M}/{m}\right)^m$ [see, for example, Devroye {\em et al.} (1996), page 218],
the inequality $1+\log x \le
2\log(1+x), \ \forall \ x\ge 1$, and taking $K_0$ such that $C_0K_0>4$ we get
\begin{eqnarray*}
\sum_{m=1}^M \binom{M}{m} \exp\left(- C_0 \frac{ m K_{1} }{a\sigma^2}
\log\left(1+\frac{M}{m}\right)
\right)
&\le&
\sum_{m=1}^M \exp\left(- m\log\left(1+\frac{M}{m}\right) \right)
\\
&\le&
\sum_{m=1}^\infty \exp (- m\log 2) < \infty.
\end{eqnarray*}
These inequalities finally yield the
bound on the tail probabilities
\begin{eqnarray}\label{eq:pd}
 &&p_\delta \le C_3 \exp\left(- C_4\frac{n\delta}{a\sigma^2}\right)
\end{eqnarray}
for some constants $C_3,C_4>0$, which easily implies the bound
(\ref{3.10.1}) on the expected value.
\end{proof}

\medskip

\begin{proof}[Proof of Theorem \ref{theorem:prima}]
Theorem \ref{theorem:prima} follows directly from Lemmas \ref{een} and
\ref{twee}.
\end{proof}

\medskip

\begin{proof}[Proof of Theorem \ref{theorem:inprob1}]

\noindent First notice that, by Lemma \ref{een}, for $a = b > 1$ there exists $C_1 > 0$  such that
\begin{eqnarray}
 && \PP \left (  \| \widetilde{f} - f \|_n^2  \geq
\inf_{\lambda\in\RR^M} \left\{ \frac{a+1}{a-1}\| {\sf f}_{\lambda} - f \|_n^2
 + \ C_1 a \sigma^2 \frac{ M({\lambda})}{n} \log\left( 1 +
\frac{M}{M({\lambda})\vee 1}\right) \right\}
 +  \delta \right ) \nonumber \\
&& \leq \PP \left ( \frac{a}{a-1}\max_{1\le m \le M} \max_{1\le k\le \binom{M}{m}} \left\{
 2a \| \Pi_{m,k}^{\lambda}W \|_n^2 - \frac{K_{1} m}{n}  \log\left(1+
\frac{M}{m\vee 1}\right) \right\} \geq \delta/2\right) \nonumber \\
&& \hspace{2cm} +  \PP \left(     \frac{2a^2}{a-1} V_n^2(\lambda)    \geq \delta/2 \right) \nonumber.
\end{eqnarray}
Next, the rescaled variable $n\sigma^{-2}V_n^2(\lambda)$ has a
$\chi^2$ distribution with 1 degree of freedom. Combining the
exponential bound for tail probabilities of $\chi^2$ random
variables (Lemma \ref{LM})  and the exponential bound (\ref{eq:pd})
completes the proof.
\end{proof}

\bigskip

\subsection{\sc Proof of Theorem \ref{th2}}
By the same reasoning as in the proof of Theorem
\ref{theorem:prima},
\begin{eqnarray*}
\| \widetilde{f} -f \|^2 &=& (1+a) \| \widetilde{f} - f  \|_n^2 + \left\{
\|\widetilde{f}-f \|^2 -(1+a) \|\widetilde{f}-f \|_n^2 \right\} \\
&\le& (1+a) \left\{ \| {\sf f}_{\lambda}-f \|_n^2
+2\left<W,\widetilde{f}-{\sf f}_{\lambda}\right>_n+\text{pen}(\lambda)-\text{pen}(\wh{\lambda})\right\} \\
&&\quad + \left\{ \|\widetilde{f}-f \|^2 -(1+a) \|\widetilde{f}-f \|_n^2 \right\}\\
&=&(1+a) \left\{ \| {\sf f}_{\lambda}-f \|_n^2
+2\left<W,\widetilde{f}-{\sf f}_{\lambda}\right>_n+\text{pen}(\lambda)-\frac{\text{pen}(\wh{\lambda})}{2}\right\} \\
&&\quad +\left\{ \|\widetilde{f}-f \|^2 -(1+a) \|\widetilde{f}-f \|_n^2 -\frac{1+a}{2}
\text{pen}(\wh{\lambda})\right\}.
\end{eqnarray*}
The first term on the right, provided $K_{1}>0$ is chosen large
enough, can be handled in exactly the same way as in the proof of
Theorem
\ref{theorem:prima}. It remains to study the second term on the right.\\

Considering separately the cases $M(\lambda)=0$ and $1\le
M(\lambda)\le M$ we obtain
\begin{eqnarray}
\nonumber
&&\|\widetilde{f}-f \|^2 -(1+a) \|\widetilde{f}-f \|_n^2 -\frac{1+a}{2}
\text{pen}(\wh{\lambda})
\\
&& \le \max \left\{U_0, \ \max_{1\le m\le M} \, \sup_{\lambda:
M(\lambda)=m}\left[ U_\lambda - \frac{1+a}{2}
\text{pen}(\lambda)\right]\right\} \nonumber
\end{eqnarray}
where $U_\lambda =\|{\sf f}_{\lambda}-f \|^2 -(1+a) \|{\sf
f}_{\lambda}-f \|_n^2 $. For each $1\le m\le M$, let the sets
$A_{m,k}(0)$, $1\le k\le \binom{M}{m}$, form a partitioning of the
set $A_{m}(0)=\{ \lambda\in \RR^M:\ M(\lambda)=m\}$.
 Deduce that, for any $\delta>0$,
\begin{eqnarray}
\label{3.11.1} && \PP\left\{ \|\widetilde{f}-f \|^2 -(1+a)
\|\widetilde{f}-f \|_n^2 -\frac{1+a}{2} \text{pen}(\wh{\lambda}) \ge
\delta \right\} \\
&&\le \PP\left\{U_0\ge \delta/2 \right\} +\sum_{m=1}^M \PP\left\{
\sup_{\lambda: M(\lambda)=m } U_\lambda \ge
D(\delta)\right\}\nonumber
\\
&&\le \PP\left\{U_0\ge \delta/2 \right\} + \sum_{m=1}^M
\sum_{k=1}^{\binom{M}{m}} \PP\left\{ \sup_{\lambda\in A_{m,k}(0)}
 U_\lambda \ge
D(\delta)\right\}\nonumber
\end{eqnarray}
where
$$
D(\delta) = \frac{(1+a)K_{1}}{2n} m \log \left( 1+ \frac{n\vee
M}{m\vee1}\right) + \frac{\delta}{2}.
$$
 The following result establishes a bound on the shatter
coefficient of the class of subgraphs of the functions $({\sf
f}_{\lambda} - f)^2$ that will be subsequently used to control the
behavior of the empirical process on the right-hand side of
(\ref{3.11.1}).

\begin{lemma}
\label{3.11} Let $\shat(n,m,k)$ be the shatter coefficient of the
collection of sets
$$\left\{ (x,\beta) :\ ({\sf f}_{\lambda}-f )^2(x) \ge \beta, \
\beta\ge0, \ x\in \X\right\},\quad \ \lambda\in A_{m,k}(0).$$ Then,
for any $1\le m\le M$, $1\le k\le \binom{M}{m}$, we have
$$\log \shat(2n,m,k)\le C m \left\{ 1+ \log\left( 1+
\frac{n}{m}\right) \right\}$$
where $C>0$ is an absolute constant.
\end{lemma}
\begin{proof}
Note that
\begin{eqnarray*}
&&\left\{ (x,\beta):\ ({\sf f}_{\lambda}-f )^2(x) \ge \beta, \
\beta\ge 0\right\} \\
&=& \left\{ (x,\beta):\ {\sf f}_{\lambda}(x)-f (x) \le - \sqrt{
\beta},  \ \beta\ge0 \right\} \cup \left \{ (x,\beta):\ {\sf
f}_{\lambda}(x) -f (x) \ge \sqrt{\beta}, \ \beta\ge0 \right\}
\end{eqnarray*}
and recall that the VC-dimension of the collection of sets $\left \{
(x,\beta):\ {\sf f}_{\lambda}(x)-f (x) \ge \sqrt{ \beta}, \
\beta\ge0 \right\} $, $\lambda\in A_{m,k}(0)$,  is less than $m+1$,
cf. Theorem 13.9 of Devroye, Gy\"orfi and Lugosi (1996)
or van de
Geer (2000), page 40.
Similarly, the VC-dimension of $\left\{ (x,\beta):\ {\sf
f}_{\lambda}(x)-f (x) \le - \sqrt{ \beta}, \ \beta\ge0 \right\} $,
$\lambda\in A_{m,k}(0)$, is less than $m+1$. Apply
 Lemma 15, page 18,  in Pollard (1984)
to deduce that the collection of sets $\left\{ (x,\beta):\ ({\sf
f}_{\lambda}-f )^2(x) \ge \beta, \ \beta\ge0 \right\}$, $\lambda\in
A_{m,k}(0)$, has VC-dimension $V_k$ less than $m+1$. The shatter
coefficient $\shat(2n,m,k)$ is related to the VC-dimension of the
latter class by the inequality
$$ \log\shat(2n,m,k) \le V_k \left\{ 1+ \log\left( 1+ \frac{2n}{V_k}\right)
\right\},$$ see, for example, Theorem 4.3 on page 145 of Vapnik
(1998). To conclude the proof, use the fact that the right-hand side
is an increasing function of $V_k$ .
\end{proof}

Now, using the inequality $D(\delta) + a\|{\sf f}_{\lambda}-f \|^2
\ge 2\sqrt{aD(\delta)}\|{\sf f}_{\lambda}-f \|$ and Theorem $5.3^*$
on page 198 of Vapnik (1998) we get
\begin{eqnarray*}
&&\PP\left\{ \sup_{\lambda\in A_{m,k}(0)}
 U_\lambda \ge
D(\delta)\right\} \\&& =\PP\left\{ \exists \lambda\in A_{m,k}(0): \
\|{\sf f}_{\lambda}-f \|\ne 0 \ \text{and} \ (1+a)\Big[\|{\sf
f}_{\lambda}-f \|^2 - \|{\sf f}_{\lambda}-f \|_n^2\Big]\ge D(\delta)
+ a\|{\sf f}_{\lambda}-f \|^2\right\}
\\
&&\le \PP\left\{ \sup_{\lambda\in A_{m,k}(0): \, \|{\sf
f}_{\lambda}-f \|\ne 0} \frac {\|{\sf f}_{\lambda}-f \|^2 -\|{\sf
f}_{\lambda}-f \|_n^2} {\|{\sf f}_{\lambda}-f \|} \ge
\frac{2\sqrt{aD(\delta)}}{1+a} \right\}\\
&&\le 4 \shat(2n,m, k)
 \exp\left\{-
 \frac{ anD(\delta) }{(1+a)^2B}\right\}.
\end{eqnarray*}
Therefore,
\begin{eqnarray*}
&&\sum_{m=1}^M \sum_{k=1}^{\binom{M}{m}} \PP\left\{ \sup_{\lambda\in
A_{m,k}(0)}
 U_\lambda \ge
D(\delta)\right\}
\\
&&\le 4\sum_{m=1}^M \sum_{k=1}^{\binom{M}{m}} \shat(2n,m, k)
 \exp\left\{-
 \frac{ anD(\delta) }{(1+a)^2B}\right\}
 \\
&&\le 4\sum_{m=1}^M {\binom{M}{m}} \exp\left\{ Cm \left[
1+\log\left( \frac{n}{m}\right)\right] \right\}
\exp\left\{- \frac{ aK_{1}m }{2(1+a)B} \log
\left( 1+ \frac{n\vee M}{m\vee1}\right) - \frac{ an\delta}{2(1+a)^2B} \right\}\\
&& \text{ by Lemma \ref{3.11} }
\\
&& \le C_5 \exp\left(-  C_6 \frac{ n\delta }{ a   B} \right) , \ \
\forall \ a>1,
\end{eqnarray*}
for $K_{1}=K_{1}(a,B)$ large enough, and some universal constants
$C_5,C_6>0$, where we have used the same crude bound for
$\binom{M}{m}$ as in the proof of Lemma \ref{twee}.
 Furthermore,
\begin{eqnarray*}
\PP\left\{
 U_0 \ge \delta/2\right\} &\le& \PP\left\{
\|f \|^2 - \|f \|_n^2\ge \frac{\sqrt{2a\delta}}{1+a}\|f \| \right\}
\\
&\le& \exp\left\{-
 \frac{ an\delta }{(1+a)^2B}\right\} \le
 \exp\left\{-
 \frac{ n\delta }{4aB}\right\}, \ \ \forall \ a>1,
\end{eqnarray*}
where the last but one inequality follows, e.g., from Proposition
2.6 in Wegkamp (2003). The exponential bounds in the last two
displays and (3.11) easily imply
\begin{eqnarray*}
\EE \left\{ \|\widetilde{f}-f \|^2 -(1+a) \|\widetilde{f}-f \|_n^2
-\frac{1+a}{2} \text{pen}(\wh{\lambda})\right\} \le C_7 \frac{Ba}{n}
\end{eqnarray*}
for some constant $C_7>0$. This concludes the proof of Theorem
\ref{th2}. \qed
\\

\bigskip

\section{Near optimal aggregation with a data dependent  $L_1$ penalty}

\noindent We consider here only the fixed design regression. In addition to Assumptions (A1) and (A2), throughout this section we suppose the
following.
\\

\noindent {\sc Assumption (A3)} {\em The matrix
 \[ \Psi_n = \left( \frac1n \sum_{i=1}^n f_j(X_i) f_{j^\prime} (X_i)
\right)_{ 1\le j,j^\prime \le M}\]
is positive definite for any given $n\ge 1$.}
\\

Let $\xi_{\min}$ be the smallest eigenvalue  of the matrix $\Psi_n$. Note that under our assumptions
\begin{equation}
\label{ksi}
0<\xi_{\min}\le \|f_j\|_n^2 \le L^2, \quad j=1,\dots,M.
\end{equation}
We propose the aggregation procedure defined by the following choice of
weights:
\begin{equation}
\label{l1pen} \wh {\lambda} = \argmin_{\lambda\in\Lambda_{M,T,2}}
\left\{ \wh{S}(\lambda)+ {\rm pen}(\lambda) \right\}
\end{equation}
where
\[ \Lambda_{M,T,2}=\left\{ \lambda\in\RR^M:\ \sum_{j=1}^{ M} \lambda_j^2
\le T^2 \right\},
\]
for $T>0$ large enough, and the penalty term is given by
\begin{eqnarray}
\label{pen5}
 \text{pen}(\lambda) =   \sum_{j=1}^M r_{n,j} |\lambda_j| \ \text{ with }
 \ \ r_{n,j} =
 2\sqrt{2}\sigma \|f_j\|_n \sqrt{\frac{ 2\log M + \log n}{n}}.
\end{eqnarray}

\begin{theorem}\label{theorem:emp}
Let $X_i \in {\mathcal X}$, $i=1,\dots,n$, be fixed. Let $\wh{\lambda}$
be the penalized least squares estimate defined by (\ref{l1pen}) with penalty
(\ref{pen5}). Set
$\widetilde{f}= {\sf f}_{\wh{\lambda}}$. Let $T>0$  be such that $T^2 \xi_{\min}> 2
L^2$.  Then, for all $a>1$, and all integers $n\ge1$, $M\ge 2$, we have,
\begin{eqnarray}\label{eq:L}
\ \qquad \EE \| \widetilde{f} - f\|_n^2 &\le&
  \inf_{\lambda\in \RR^M}
\left\{ \frac{a+1}{a-1} \| {\sf f}_{\lambda} - f\|_n^2  + \frac{16a^2}{a-1}
\left(\frac{ \sigma^2 L^2}{\xi_{\min}} \right) \frac{2\log M+\log n}{n}
M(\lambda)
\right\}
\\
&& \qquad + \ \frac{(T+M^{-1/2})^2L^2}{n\sqrt{\pi(2\log M+\log n)}}.
\nonumber
\end{eqnarray}
\end{theorem}

\medskip

\begin{corollary}[MS]\label{cor40}
Let assumptions of Theorem \ref{theorem:emp} be satisfied and $T\le
(\log(M\vee n))^{1/4}$. Then there exists a constant
$C=C(T,L,\sigma^2,\xi_{\min})>0$ such that for all $\eps>0$ and for
all integers $n\ge 1$ and $M\ge 2$,
\begin{eqnarray}
\EE \| \widetilde{f} - f \|_n^2   \le (1+\eps)\inf_{1\le j\le M} \|
f_j - f \|_n^2 + C \left(1+\eps+{\eps}^{-1}\right) \frac{ \log
(M\vee n)} {n}. \nonumber
\end{eqnarray}
\end{corollary}
\begin{proof}[Proof]
Using assumptions on $T$ and (\ref{ksi}), we trivially get $T
>\sqrt{2 L^2/\xi_{\min}} \ge M^{-1/2}$. This implies that the last
summand in (\ref{eq:L}) is $O(1/n)$. The rest of the proof is
analogous to that of Corollary \ref{cor:MS}.
\end{proof}
\medskip

\begin{corollary}[C]\label{cor41}
Let assumptions of Theorem \ref{theorem:emp} be satisfied and $T\le
(\log(M\vee n))^{1/4}$. Then there exists a constant
$C=C(T,L,\sigma^2,\xi_{\min})>0$  such that for all $\eps>0$ and for
all integers $n\ge 1$ and $M\ge 2$,
\begin{eqnarray}
\EE \| \widetilde{f} - f \|_n^2 \le
(1+\eps)\inf_{\lambda\in\Lambda^M} \| {\sf f}_{\lambda} - f \|_n^2 +
C \left(1+\eps +{\eps}^{-1}\right) \overline{\psi}_n^C(M), \nonumber
\end{eqnarray}
where
\begin{eqnarray*} \overline{\psi}_n^C(M) = \begin{cases} (M\log n)/n &
\text{ if }
M\le\sqrt{n},\\
\sqrt{(\log M ) / n}   & \text{ if } M> \sqrt{n}.
\end{cases}\end{eqnarray*}
\end{corollary}
\begin{proof}[Proof]
We bound the last summand in (\ref{eq:L}) as in the previous proof and
we use then the argument similar to that of the
proof of Corollary \ref{cor:C}.
\end{proof}
\medskip

\begin{corollary}[L]\label{cor42}
Let assumptions of Theorem \ref{theorem:emp} be satisfied and $T\le
(\log(M\vee n))^{1/4}$. Then there exists a constant
$C=C(T,L,\sigma^2,\xi_{\min})>0$ such that for all $\eps>0$ and for
all integers $n\ge 1$ and $M\ge 2$,
\begin{eqnarray}
\EE \| \widetilde{f} - f \|_n^2 \le (1+\eps)\inf_{\lambda\in\RR^M}
\| {\sf f}_{\lambda} - f  \|_n^2 + C \left(1+\eps
+{\eps}^{-1}\right)\frac{ M \log (M\vee n)}{n}. \nonumber
\end{eqnarray}
\end{corollary}
\begin{proof}[Proof]
We bound the last summand in (\ref{eq:L}) as in the proof Corollary
\ref{cor40} and we use that $M(\lambda)\le M$.
\end{proof}
\medskip

\begin{proof}[Proof of Theorem \ref{theorem:emp}]
 We begin as in  Loubes and Van de Geer (2002).
By definition, $\widetilde{f}= {\sf f}_{\widehat\lambda}$ satisfies
\begin{eqnarray*}
\wh{S}(\widehat\lambda) +  \sum_{j=1}^M r_{n,j} |\widehat \lambda_j|  \le
\wh{S}(\lambda) + \sum_{j=1}^M r_{n,j} |\lambda_j|
\end{eqnarray*} for all $\lambda\in\Lambda_{M,T,2}$,
which we may rewrite as
\begin{eqnarray*}
\| \widetilde{f} - f\|_n^2 +  \sum_{j=1}^M r_{n,j} |\widehat \lambda_j| \le
\| {\sf f}_{\lambda} - f\|_n^2 +  \sum_{j=1}^M r_{n,j} |\lambda_j| +2\left< W,
\widetilde{f} - {\sf f}_{\lambda} \right>_n.
\end{eqnarray*}
We define the random variables \[  V_j =  \frac1n \sum_{i=1}^n f_j(X_i)
W_i,\quad 1\le j\le M,\] and the event
\[ A = \bigcap_{j=1}^M \left\{  2| V_j | \leq r_{n,j} \right\}.\]
The normality assumption (A1) on $W_i$ implies that $\sqrt{n}\, V_j\sim
N\left(0,\sigma^2\|f_j\|_n^2 \right)$, $1\le j\le M$. Applying the union
bound followed by the standard tail bound for the $N(0,1)$ distribution,
yields
\begin{eqnarray}\label{eq:pac}
\PP ( A^c)
&\le&
\sum_{j=1}^M  \PP\{ \sqrt{n}|V_j| > \sqrt{n}r_{n,j}/2\}
\le
\sum_{j=1}^M\frac{4 }{\sqrt{2\pi}} \frac{\sigma \|f_j\|_n}{ \sqrt{n}
r_{n,j}} \exp\left(
-\frac{n r_{n,j}^2}{ 8 \sigma^2 \|f_j\|_n^2} \right)\\
&=&
\frac{1}{Mn \sqrt{\pi(2\log M+\log n)}}. \nonumber
 \end{eqnarray}
Then, on the set $A$, we find
\begin{eqnarray*}
2 \left<W, \widetilde{f} - f \right>_n = 2\sum_{j=1}^M V_j (\widehat
\lambda_j - \lambda_j) \le   \sum_{j=1}^M r_{n,j} |\widehat  \lambda_j -
\lambda_j|
\end{eqnarray*}
and therefore, still on the set $A$,
\begin{eqnarray*}
\| \widetilde{f} - f\|_n^2 & \le&  \| {\sf f}_{\lambda} - f\|_n^2  +  \sum_{j=1}^M
r_{n,j}|\widehat \lambda_j - \lambda_j| +  \sum_{j=1}^M r_{n,j}
|\lambda_j| -  \sum_{j=1}^M r_{n,j} |\widehat\lambda_j|.
\end{eqnarray*}
Recall that $J(\lambda)$ denotes the set of indices of the non-zero
elements of $\lambda$, and   $M(\lambda) = {\rm Card\ } J(\lambda)$.
Rewriting the right-hand side of the previous display, we find, on the
set $A$,
\begin{eqnarray*}
\| \widetilde{f} - f\|_n^2 & \le&  \| {\sf f}_{\lambda} - f\|_n^2  + \left(
\sum_{j=1}^M r_{n,j} |\widehat \lambda_j - \lambda_j| -  \sum_{j \not\in
J(\lambda)} r_{n,j} |\widehat\lambda_j| \right) \\&&\qquad\qquad\quad +
\left(  -  \sum_{j \in J(\lambda)} r_{n,j}|\widehat\lambda_j|
 +  \sum_{j\in J(\lambda)} r_{n,j}|\lambda_j| \right)
 \\
&\le&  \| {\sf f}_{\lambda} - f\|_n^2  +2 \sum_{j\in J(\lambda) }
r_{n,j}|\widehat \lambda_j - \lambda_j|
\end{eqnarray*}
by the triangle inequality and the fact that $\lambda_j=0$ for $j\not\in
J(\lambda)$. Since $\xi_{\min}>0$, we have
\begin{eqnarray*}
\xi_{\min}^{-1}\| \widetilde{f} - {\sf f}_{\lambda}\|_n^2 \ge  \sum_{j\in J(\lambda)
} |\widehat \lambda_j - \lambda_j|^2.
\end{eqnarray*}
Combining this with the Cauchy-Schwarz and triangle inequalities,
respectively,
we find further that, on the set $A$,
\begin{eqnarray}\label{eq:bo}
\label{ontheset}
\| \widetilde{f} - f\|_n^2 & \le&   \| {\sf f}_{\lambda} - f\|_n ^2  + 2 \sum_{j\in
J(\lambda) } r_{n,j}|\widehat \lambda_j - \lambda_j| \\
& \le&  \| {\sf f}_{\lambda} - f\|_n ^2  +
2\sqrt{\xi_{\min}^{-1}}\sqrt{\sum_{j\in J(\lambda)} r_{n,j}^2} \left( \|
\widetilde{f} -
f\|_n + \| {\sf f}_{\lambda} -f \|_n \right)\nonumber\\
&\le&
\nonumber
\| {\sf f}_{\lambda} - f\|_n ^2  + 2\sqrt{\xi_{\min}^{-1}} r_n \sqrt{ M(\lambda)
} \left( \| \widetilde{f} - f\|_n + \| {\sf f}_{\lambda} -f \|_n \right),
\end{eqnarray}where
\[ r_n \egal 2\sqrt{2} \, L\sigma  \sqrt{\frac{2\log M+\log n}{n}}.\]
Inequality (\ref{ontheset}) is of the simple form $ v^2 \le c^2 + vb + cb $
with $v=\|\widetilde{f}-f\|_n$, $b=2 r_n \sqrt{M(\lambda)/\xi_{\min}}$ and
$c=\|{\sf f}_{\lambda}-f\|_n$.
 After applying the inequality $2xy\le x^2/\alpha+ \alpha y^2$
($x,y\in\RR, \ \alpha>0$) twice, to $2bc$ and $2bv$, respectively, we easily find
 $v^2 \le  v^2/(2\alpha) + \alpha\, b^2 +
(2\alpha+1)/(2\alpha)\, c^2 $, whence $v^2\le a/(a-1) \{ b^2 (a/2) + c^2
(a+1)/a \}$
for $a=2\alpha>1$. Recalling that (\ref{ontheset}) is valid on the set
$A$, we now get
 that
\[ \EE \left[\| \widetilde{f} -f\|_n ^2 I_A\right]
\le \inf_{\lambda\in\Lambda_{M,T,2}} \left\{ \frac{a+1}{a-1} \| {\sf f}_{\lambda}
- f\|_n^2 +
\frac{2a^2}{\xi_{\min}(a-1)}
r_n^2 M(\lambda) \right\} , \quad \forall \ a>1.
\]
Consequently, since by the Cauchy-Schwarz inequality,
$$\| \widetilde{f} - f\|_\infty\le L(\sum_{j=1}^M |\lambda_j|+ 1)
\le( \sqrt{M}T+1)L,$$
 we find
\begin{eqnarray}\label{eq:inf} \EE \| \widetilde{f} -f\|_n ^2
&\le& \EE\left[\| \widetilde{f}-f\|_n^2 I_A\right] + (\sqrt{M} T+1)^2L^2 \PP(A^c)
\nonumber\\
&\le& \inf_{\lambda\in\Lambda_{M,T,2}} \left\{ \frac{a+1}{a-1} \|
{\sf f}_{\lambda} - f \|_n^2 +
\frac{2a^2 r_n^2 }{(a-1)\xi_{\min}}
M(\lambda) \right\}\\ && \quad +
\ \frac{(T+M^{-1/2})^2L^2}{n\sqrt{\pi(2\log M+\log n)}}.
\nonumber
\end{eqnarray}
It remains to show that (\ref{eq:inf}) remains valid with the set
$\Lambda_{M,T,2}$
replaced by the entire $\RR^M$. For this, observe that
$\lambda\not\in \Lambda_{M,T,2}$ implies $\sum_{j=1}^M\lambda_j^2 > T^2$,
and thus $\| {\sf f}_{\lambda}\|_n^2 \ge \xi_{\min}\sum_{j=1}^M\lambda_j^2
> \xi_{\min} T^2$. Therefore, for $\lambda\not\in \Lambda_{M,T,2}$, we
have $$\| {\sf f}_{\lambda}-f\|_n \ge \| {\sf f}_{\lambda}\|_n  -\|f\|_n >
\sqrt{\xi_{\min}} T  - L> L$$ by our choice of $T$. On the other hand,
for  $\lambda=0 \in \Lambda_{M,T,2}$, we have $$\| {\sf f}_{\lambda}-f\|_n = \|
f\|_n \le L$$ and  pen$(0)=0$. Thus, the value of the whole expression
under the infimum in (\ref{eq:inf}) for $\lambda=0$ is strictly smaller
than the value of this expression for any $\lambda\not\in
\Lambda_{M,T,2}$, which proves the result.
\end{proof}

As in Section 3.1, we present now a statement in probability that
complements the results of this section.

\begin{theorem}\label{theorem:empprob}
Let $X_i \in {\mathcal X}$, $i=1,\dots,n$, be fixed. Let $\wh{\lambda}$
be the penalized least squares estimate defined by (\ref{l1pen}) with $ \Lambda_{M,T,2}$ replaced by $\RR^M$ and  with penalty
(\ref{pen5}). Set
$\widetilde{f}= {\sf f}_{\wh{\lambda}}$. Then, for all $a>1$, and all integers $n\ge1$, $M\ge 2$, we have,
\begin{eqnarray}\label{eq:L1}
&& \PP\left (  \| \widetilde{f} - f\|_n^2 \geq
  \inf_{\lambda\in \RR^M}
\left\{ \frac{a+1}{a-1} \| {\sf f}_{\lambda} - f\|_n^2  + \frac{16a^2}{a-1}
\left(\frac{ \sigma^2 L^2}{\xi_{\min}} \right) \frac{2\log M+\log n}{n}
M(\lambda)
\right\} \right)
\\
&& \qquad \leq  \frac{1}{Mn \sqrt{\pi(2\log M+\log n)}}.
\nonumber
\end{eqnarray}
\end{theorem}

\begin{proof}

\noindent This result follows directly from the proof of Theorem
\ref{theorem:emp}. Note first that now (\ref{eq:bo}) is valid for
all $\lambda \in \RR^M$ and not only for $\lambda
\in\Lambda_{M,T,2}$. Using (\ref{eq:bo}) and the argument after it
we find that the left hand side in (\ref{eq:L1}) can be bounded by
$\PP(A^c)$. The result follows by invoking  (\ref{eq:pac}).

\end{proof}

\medskip

\noindent {\sc Remarks.}
\\

\noindent 1. The method  presented in this section is not strictly
an $L_1$-penalized one. Indeed, it implements two penalties: the
data dependent $L_1$-penalty $\sum_{j=1}^M r_{n,j} |\lambda_j|$, and
the $L_2$-penalty $\sum_{j=1}^M\lambda_j^2$ that appears implicitly
via the choice of the set $\Lambda_{M,T,2}$. The resulting
minimization problem can be solved in practice using standard convex
programming software. The $L_2$ part of the penalty is less
influential, since it should  typically be applied with $T\to\infty$
as $M$ (respectively $n$) grows, which means that the restriction to
$\Lambda_{M,T,2}$ becomes asymptotically negligible. Moreover, the
restriction is not always needed. For example, the bound in
probability (Theorem \ref{theorem:empprob}) is obtained for
$\wh{\lambda}$ that minimizes the $L_1$-penalized least squares over
the entire $\RR^M$.
\\

\noindent 2. Assumption (A3) is mild,  and it is also made by Efron {\em et al.} (2004) in the context of LARS. In practice, this assumption can
 always be  checked. A stronger assumption is that $\xi_{\min}>c$ for some constant $c>0$, independent of $n$ and $M$ if one or both of these parameters
are allowed to grow (which is typically the more interesting case). There are at least two important examples where such a stronger assumption
holds. The first example is standard in the parametric regression context: $M$ is fixed and $\Psi_n/n \to \Psi$ where $\Psi$ is a nonsingular
$M\times M$ matrix. The second one is related to nonparametric regression: $M=M_n$ is allowed to go to $\infty$ as $n\to\infty$ and the
functions $f_j$ are orthogonal with respect to the empirical norm. This corresponds, for instance, to sequence space models, where the
estimators $f_j=\wh{f_j}$ are constructed from non-intersecting blocks of coefficients. Aggregating such mutually orthogonal estimators may lead
to adaptive estimators with good asymptotic properties [{\em cf., e.g.}, Nemirovski (2000)]. Local image smoothing provides us an application
where the condition $\xi_{\min}>c$ is naturally satisfied. For example, Katkovnik {\em et al.} (2002, 2004) suggest different methods of
aggregation of local image estimators obtained from non-intersecting sectors around a given pixel (these estimators are mutually orthogonal with
respect to the
empirical norm).\\

\noindent 3. Inspection of the proofs shows that the constants $C=C(T,L,\sigma^2,\xi_{\min})$ in Corollaries \ref{cor40}, \ref{cor41},
\ref{cor42} have the form $C= A_1 + A_2\xi_{\min}^{-1}$, where $A_1$ and $A_2$ are constants independent of $\xi_{\min}$. In general,
$\xi_{\min}$ may depend on $n$ and $M$. However, if $\xi_{\min}>c$ for some constant $c>0$, independent of $n$ and $M$, as previously discussed,
the rates of aggregation given in Corollaries \ref{cor40}, \ref{cor41}, \ref{cor42} are near optimal, up to logarithmic factors. They are even
exactly optimal ({\em cf.} (\ref{oragg}) and the lower bounds of the next section) for some configurations of $n,M$: for (MS)-aggregation if
$n^{a'} \le M\le n^{a}$, and for (C)-aggregation if $n^{1/2} \le M\le n^{a}$, where $0<a'<a<\infty$.
\\

\noindent 4. From the bound in Theorem \ref{theorem:emp}, we see that $T$ is allowed to grow with $n$ and $M$ (as fast as $T\asymp (\log(M\vee
n))^{1/4}$ is possible). Moreover, the proof of Theorem \ref{theorem:emp} reveals that by taking a larger constant than $2\sqrt{2}$ in
(\ref{pen5}), even faster rates are allowed, for example, $T$ can grow as a power of $n$. This may be needed to guarantee the condition $T^2 > 2
L^2 / \xi_{\min}$ for $n$ large enough, because the value $L$ is typically not known and $\xi_{\min}$ may depend on $n$ and $M$. However, the
condition $T^2
> 2 L^2 / \xi_{\min}$ is only needed to cover the linear
aggregation. For (MS) and (C) aggregation, Corollaries \ref{cor40}, \ref{cor41} can be obtained directly from (\ref{eq:inf}), and thus it
suffices  to take any $T\ge 1$, since $\Lambda^M \subset \Lambda_{M,1,2}$, or to replace $\Lambda_{M,T,2}$ by $\Lambda^M$ in the definition of
$\wh{\lambda}$.

\section{Lower bounds}

\noindent For regression with random design and the $L_2(\RR^d,
\mu)$-risks, lower bounds for aggregation and optimal rates
$\psi_{n,M}$ as given in (\ref{oragg}) were established by Tsybakov
(2003). In this section we extend the lower bounds of Tsybakov
(2003) for (MS) and (L) aggregation to regression with fixed design.
Further, we state these bounds in a more general form, considering
not only the expected squared risks, but also other loss functions.
This generalization allows one to treat optimality of the upper
bounds ``in probability" obtained in the previous sections (Theorems
\ref{theorem:inprob1}, \ref{theorem:empprob}). It shows that the
remainder terms in these bounds are optimal or near optimal for the
(MS) and (L) aggregation.

In this section we suppose that $X_1,\dots, X_n$ are fixed and that $M\le n$. Let $w:\RR \to [0,\infty)$ be a {\it loss function}, {\em i.e.}, a
monotone non-decreasing function satisfying $w(0)=0$ and $w\not \equiv 0$.

\begin{theorem}
\label{theorem:lower} Let $X_i \in {\mathcal X}$, $i=1,\dots,n$, be
fixed and $2\le M\le n$. Assume that $H^M$ is either the whole
$\RR^M$ (the (L) aggregation case) or the set of vertices of
$\Lambda^M$ (the (MS) aggregation case). Let the corresponding
$\psi_{n,M}$ be given by (\ref{oragg}) and let $M\log M\le n$ for
the case of (MS) aggregation. Then there exist $f_1,\dots,f_M\in
{\mathcal F}_0$ such that, for any loss function $w(\cdot)$,
\begin{equation}
\label{lb}
\inf_{T_n}\sup_{f\in {\mathcal F}_0}
\EE w
\Big[ \psi_{n,M}^{-1}
\Big(\|T_n - f\|^2_n - \inf_{\lambda\in H^M} \|{\sf f}_\lambda - f\|^2_n\Big)
\Big]
\ge c,
\end{equation}
where $\inf_{T_n}$ denotes the infimum over all estimators and
the constant $c>0$ does not depend on $M$ and $n$.
\end{theorem}

\medskip

Setting $w(u) = u$ in Theorem \ref{theorem:lower} we get the lower bounds for expected squared risks showing optimality or near optimality of
the remainder terms in the oracle inequalities of Corollaries \ref{cor:MS}, \ref{corL},  \ref{cor40},  \ref{cor42}. The choice of $w(u)=
I\{u>a\}$ with some fixed $a>0$ leads to the lower bounds for probabilities showing near optimality of the remainder terms in the corresponding
upper bounds (see Theorems \ref{theorem:inprob1}, \ref{theorem:empprob}).

\begin{proof}
We proceed similarly to Tsybakov (2003). The proof is based on the following lemma [which can be obtained, for example, by combining Theorems
2.2 and 2.5 in Tsybakov (2004)].

\begin{lemma}
\label{l1} Let $w$ be a loss function, $A>0$ be such that $w(A)>0$, and
let ${\mathcal C}$ be a finite set of functions on ${\mathcal X}$ such that $N={\rm card}({\mathcal C})\ge 2$,
$$\|f-g\|_n^2 \ge 4s^2 >0, \quad \forall \ f,g \in {\mathcal C}, \ \ f\neq g,$$
and the Kullback divergences $K(\PP_{f}, \PP_{g})$ between the measures
$\PP_{f}$ and $\PP_{g}$ satisfy
$$
K(\PP_{f}, \PP_{g}) \le (1/16) \log N, \quad \forall \ f,g \in {\mathcal C}.
$$
Then for $\psi=s^2/A$ we have
$$
\inf_{T_n} \sup_{f\in {\mathcal C}}
\EE w
\Big[ \psi^{-1}
\|T_n - f\|^2_n
\Big]
\ge c_1 w(A),
$$
where $\inf_{T_n}$ denotes the infimum over all estimators and $c_1>0$ is a constant.
\end{lemma}

\noindent {\sl The (MS) aggregation case.} Let $H^M$ be the set of
vertices of $\Lambda^M$, $M\log M\le n$, and $\psi_{n,M}= (\log
M)/n$. Pick $M$ disjoint subsets $S_1,\dots, S_M$ of $\{X_1,\dots,
X_n\}$, each $S_j$ of cardinality $\log M$ (w.l.o.g. we assume that
$\log M$ is an integer) and define the functions
$$
f_j(x) = \gamma I\{x\in S_j\}, \quad j=1,\dots,M,
$$
where $\gamma\le L$ is a positive constant to be chosen. Clearly, $\{f_1,\dots, f_M\} \subset {\mathcal F}_0$. Thus, it suffices to prove the
lower bound of the theorem where the supremum over $f\in {\mathcal F}_0$ is replaced by that over $f \in \{f_1,\dots, f_M\}$. But for such $f$
we have $\min_{1\le j\le M} \|f_j - f\|^2_n=0$, and to finish the proof for the (MS) case, it is sufficient to bound from below the quantity $
\sup_{f \in \{f_1,\dots, f_M\}} \EE w ( \psi_{n,M}^{-1} \|T_n - f\|^2_n )$, where $\psi_{n,M}= (\log M)/n$, uniformly over all estimators $T_n$.
This is done by applying Lemma \ref{l1}. In fact, note that, for $j\neq k$,
\begin{equation}
\label{3}
\|f_j-f_k\|^2_n = \frac{2\gamma^2\log M}{n} \egal 4s^2.
\end{equation}
Since
$W_j$'s are ${
N}(0,\sigma^2)$ random variables, the Kullback divergence
$K(\PP_{f_j}, \PP_{f_k})$ between $\PP_{f_j}$ and $\PP_{f_k}$
satisfies
\begin{equation}
\label{kull}
K(\PP_{f_j}, \PP_{f_k}) = {n \over 2\sigma^2} \|f_j-f_k\|_n^2, \quad j=1,\dots,M.
\end{equation}
 In view of (\ref{3}) and
(\ref{kull}), one can choose $\gamma$ small enough to have $K(\PP_{f_j}, \PP_{f_k})\le (1/16)\log M$ for $j,k =1,\dots,M$. Now, to get the lower
bound for the (MS) case, it remains to use this inequality, identity (\ref{3}) and Lemma \ref{l1}.\\

\noindent {\sl The (L) aggregation case.} Let $H^M=\RR^M$ and $\psi_{n,M}= M/n$. Define the functions $f_j=\gamma I\{x=X_j\}, \quad
j=1,\dots,M,$ with $0<\gamma\le L$ and introduce a finite set of their linear combinations
\begin{equation}
\label{10bis} {\mathcal U} = \Big\{g=\sum_{j=1}^M \omega_j f_j: \omega\in\Omega\Big\},
\end{equation}
where
$\Omega$ is the set of all vectors $\omega\in \RR^M$ with binary coordinates
$\omega_j\in\{0,1\}$. Since the supports of $f_j$'s are disjoint, the functions
$g\in {\mathcal U}$ are uniformly bounded by $\gamma$, thus ${\mathcal U}\subset {\mathcal F}_0$.
Clearly,
$\min_{\lambda\in\RR^M} \|{\sf f}_\lambda - f\|^2_n =0$ for any $f\in {\mathcal U}$.
Therefore, similarly to the (MS) case,
it is
sufficient to bound
from below the quantity $ \sup_{f \in {\mathcal U}}
\EE w
( \psi_{n,M}^{-1} \|T_n - f\|^2_n )$ where $\psi_{n,M}= M/n$,
uniformly over all estimators $T_n$.

Note that for any $g_1=\sum_{j=1}^M \omega_j f_j \in {\mathcal U}$ and
$g_2=\sum_{j=1}^M \omega_j' f_j \in {\mathcal U}$ we have
\begin{equation}
\label{9}
\|g_1-g_2\|^2_n = \frac{\gamma^2}{n} \sum_{j=1}^M (\omega_j - \omega_j')^2 \le \gamma^2M/n.
\end{equation}
Let first $M\ge 8$. Then it follows from the Varshamov-Gilbert bound
(see, for instance, Tsybakov (2004), Chapter 2) that there exists a
subset ${\mathcal U}_0$ of ${\mathcal U}$ such that ${\rm
card}({\mathcal U }_0)\ge 2^{M/8}$ and
\begin{equation}
\label{10}
\|g_1-g_2\|^2_n \ge C_1\gamma^2M/n.
\end{equation}
for any $g_1,g_2 \in {\mathcal U}_0$. Using (\ref{kull}) and
(\ref{9}) we get, for any $g_1,g_2 \in {\mathcal U}_0$,
$$
K(\PP_{g_1}, \PP_{g_2}) \le C_2\gamma^2 M \le C_3 \gamma^2\log ({\rm
card}({\mathcal U}_0)),
$$
and by choosing $\gamma$ small enough, we can finish the proof in
the same way as in the (MS) case. If $2\le M\le 8$, we have
$\psi_{n,M}\le 8/n$, and the proof is easily obtained by choosing
$f_1\equiv 0$ and $f_2\equiv \gamma n^{-1/2}$ and applying Lemma
\ref{l1} to the set ${\mathcal U}_0=\{f_1,f_2\}$.
\end{proof}

\bigskip

\appendix

\section{}

\begin{lemma}
\label{l2} Let $f,f_1,\dots,f_M\in {\mathcal F}_0$ and $1\le m\le
M$. Let $\C$ be the finite set of functions defined in the proof of
Corollary \ref{cor:C}. Then (\ref{one}) holds and
\begin{eqnarray}
\label{one0} \min_{g\in\C} \| g-f\|^2 \le \min_{\lambda\in
\Lambda^M} \| {\sf f}_{\lambda}- f\|^2 + \frac{L^2}{m}.
\end{eqnarray}
\end{lemma}
\begin{proof}
Let $f^*$ be the minimizer of $\| {\sf f}_\lambda-f\|^2$ over
$\lambda\in\Lambda^M$. Clearly, $f^*$ is of the form
\[ f^*=\sum_{j=1}^M p_j f_j \ \text{ with } p_j\ge 0 \ \text{ and }
\sum_{j=1}^M p_j \le 1.\] Define a probability distribution on
$j=0,1,\dots,M$ by
\begin{eqnarray*} \pi_j = \begin{cases} p_j &
\text{ if }
j\neq 0,\\
1- \sum_{j=1}^M p_j  & \text{ if } j=0.
\end{cases}\end{eqnarray*}
Consider $m$ i.i.d. random integers $j_1,\ldots,j_m$ where each
$j_k$ is distributed according to $\{\pi_j\}$ on $\{0,1,\dots,M\}$.
Introduce the random function
\[
{\bar f}_m=\frac1m \sum_{k=1}^m g_{j_k}
\]
where
\begin{eqnarray*} g_j = \begin{cases} f_j &
\text{ if }
j\neq 0,\\
0  & \text{ if } j=0.
\end{cases}\end{eqnarray*}
For every $x\in {\mathcal X}$ the random variables
$g_{j_1}(x),\dots,g_{j_m}(x)$ are i.i.d. with
$\mathbb{E}(g_{j_k}(x))=f^*(x)$. Thus,
\begin{eqnarray*}
\mathbb{E} ( {\bar f}_m(x)-f^*(x))^2 &=&
 \mathbb{E}\left(\left[\frac1m \sum_{k=1}^m \{g_{j_k}(x)-
 \mathbb{E}(g_{j_k}(x))\}\right]^2 \right)\\ &\le & \frac1m
 \mathbb{E}(g_{j_1}^2(x))
\le \frac{L^2}{m} \ .
\end{eqnarray*}
Hence for every $x\in {\mathcal X}$ and every $f\in {\mathcal F}_0$
we get
\begin{eqnarray}\label{maur}
\mathbb{E} ( {\bar f}_m(x)-f(x))^2 &=& \mathbb{E} ( {\bar
f}_m(x)-f^*(x))^2 + ( f^*(x)-f(x))^2  \\
&\le & \frac{L^2}{m} + ( f^*(x)-f(x))^2 . \nonumber
\end{eqnarray}
Integrating (\ref{maur}) over $\mu(dx)$ and recalling the definition
of $f^*$ we obtain
\begin{eqnarray}\label{maur1}
\mathbb{E}\| {\bar f}_m-f\|^2 &\le& \min_{\lambda\in\Lambda^M}
\|f_\lambda-f\|^2 +\frac{L^2}{m}.
\end{eqnarray}
Finally, note that the random function ${\bar f}_m$ takes its values
in $\C$, which implies that
$$
\mathbb{E}\| {\bar f}_m-f\|^2 \ge \min_{g \in \C} \| g - f\|^2.
$$
This and (\ref{maur1}) prove (\ref{one0}). The proof of (\ref{one})
is analogous, with the only difference that (\ref{maur}) is
integrated over the empirical measure rather than over $\mu(dx)$.
\end{proof}

\end{document}